# Mathematical Modelling of Hypothalamus-Pituitary-Adrenal Axis Dynamics: A Review, a Novel Approach and Future Directions

Imran Ahmad[1], Itishree Jena[2], A. Priyadarshi*[2]

DST-CIMS, Banaras Hindu University, Varanasi-221005

[1]Department of Mathematics, Institute of Science

[2]Banaras Hindu University, Varanasi-221005 India

*Email: anupampriya@bhu.ac.in

**Abstract:**

A novel mathematical model for the hypothalamic–pituitary–adrenal (HPA) axis is proposed to comprehend the oscillations observed in hormone concentration and potential dysfunction within the HPA axis in stressful situation. This model integrates impact of hippocampal receptors on the secretion of corticotropin-releasing hormone (CRH), an additional signalling pathway involving Arginine Vasopressin (AVP) for the production and secretion of adrenocorticotropic hormone (ACTH), the inclusion of a daylight-related function for modelling circadian rhythms, and a short negative feedback loop from the pituitary to the hypothalamus in a minimal mechanistic model of the HPA axis. This expansion allows us to estimate model parameters that led to a significant reduction in the mean absolute percent error, thereby enhancing the model's predictive accuracy with a demonstration of a strong fit to the validation dataset. Through sensitivity and correlation analyses, our study shows the parameters that exert the most significant influence on the dynamics of cortisol within the system. This study reveals intricate interdependencies within the model and among the various components of the HPA axis. Despite current advancements in the model to comprehend the dynamics of the HPA Axis, these models are still lacking in addressing major factors such as the impact of genetic and epigenetics, the role of the amygdala in processing and transmitting stress signals, and other minor impacting factors.

**Glossary:**

    SCN -Suprachiasmatic nucleus

    PVN -Paraventricular nucleus

    CRH -Corticotrophin-releasing hormone

    ACTH –Adrenocorticotropin

    AVP -Arginine vasopressin

    MR -Mineralocorticoid receptors

    GR -Glucocorticoid receptors

    PTSD -Post-traumatic stress disorder

    CBG -Cortisol-binding globulin

    RMSE -Root mean square error

    MAPE -Mean absolute percent error

**1. Introduction**

The present world is full of chaos where combination of physiological and psychological factors cause continuous stress which leads depression in human beings and ultimately may lead to death in cardiovascular disorders including heart attack and heart failure. These diseases stem from the activation *hypothalamic-pituitary-adrenal (HPA)* axis. The dysregulation of HPA axis in response to chronic stress often disrupt on mental and emotional well-being by causing deep concern of *Post-traumatic Stress Disorder (PTSD)*, major depressive disorder, alcoholism, and anorexia nervosa [1-3]. To address these social psychological challenges, Mathematical HPA axis model hold potential for innovative, evidence-based solutions.

**1.1 Physiology of HPA Axis**

The HPA axis is part of the neuro-endocrine system responsible for regulating cortisol hormone. This regulation occurs through a combination of feedforward and feedback mechanisms to control various bodily processes in both normal and stressful conditions by managing the levels of corticosteroids produced by the adrenal glands. The HPA axis can be activated in three ways: (i) activation is due to the internal biological clock located in the hypothalamic suprachiasmatic nucleus (SCN) which maintains a rhythmic activity in the HPA system, creating a circadian rhythm [4–6]. (ii) the pulsatile release of hormones and the negative feedback loops within the HPA axis contribute to its activation [7]. (iii) stressful situations can also activate the HPA axis, leading to an additional release of cortisol, superimposed on the circadian and ultradian rhythms [8]. Stress-related physiological signals, such as inflammation or low blood pressure, trigger the activation of neurons located in the *paraventricular nucleus* (PVN) of the hypothalamus. Later on, these neurons release a variety of hormones, including corticotrophin-releasing hormone (CRH) and arginine vasopressin (AVP), into the portal circulation. Consequently, this leads to the stimulation of adrenocorticotropin (ACTH) release from corticotroph cells situated in the anterior pituitary gland [9]. It is further transported by blood circulation, especially to the adrenal cortex, where it stimulates the production and release of corticosteroids (cortisol in humans, and corticosterone in animals) [7, 8]. Once cortisol enters the bloodstream, it plays a crucial role in regulating negative feedback to maintain the body's balance during times of stress. However, when this feedback system is disrupted, it can lead to a decrease in the sensitivity of the HPA axis [10]. This happens because cortisol has a stronger affinity for mineralocorticoid receptors (MR) in comparison to glucocorticoid receptors (GR). When cortisol binds to MR, it forms a complex, which then undergoes homodimerization to increase its activity [11]. In particular, the GR-cortisol complex binds to CRH in the hypothalamus and ACTH in the pituitary gland to inhibit the further production of these hormones, creating a negative feedback loop cycle [12, 13]. In addition to it, a short feedback loop exists to cause ACTH to inhibit the release of CRH from the posterior pituitary gland. The hippocampus plays a critical role in the functioning of the HPA axis by stimulating the hypothalamus to produce CRH in response to stressful situations, with cortisol providing feedback on MR and GR receptors in the hippocampus [14-16]. It is well known that amount of cortisol binding to MR in comparison to GR receptors can determine whether the involvement of hippocampal mechanisms results in an overall positive or negative feedback mechanism. Study [17] revealed that hippocampal GR receptors contribute to positive feedback in the HPA axis, while hippocampal MR receptors contribute to negative feedback. It is well known that stress has a tendency to disrupt the balance of cortisol in form of oscillation in its concentration. Although an effective HPA axis has an ability to restore hormone levels to their baseline after a certain period, it is important to maintain cortisol within a specific range, as elevated levels (hypercortisolism) can give rise to conditions like depression, diabetes, visceral obesity, and osteoporosis [14, 18]. Also, insufficient cortisol levels can disturb memory formation or even lead to life-threatening adrenal crisis [19, 20]. The regulation of hormones by the HPA axis plays a pivotal role in preserving one's health.

Circadian rhythm is influenced by external factors such as daylight, temperature, and both psychological and physical stress. While, ultradian rhythm is the inherent dynamics of the HPA axis [21-25]. Stress tends to increase the amplitude of these oscillations, while frequency is unaltered by stress [26]. A human being is a diurnal creature, exhibits an autonomous circadian activity of the HPA axis that is synchronized with the 24-hour light/dark cycle. Cortisol levels in human are lowest during sleep between 8 p.m. and 2 a.m. and gradually rise as approach to peak in morning hour between 6 a.m. and 10 a.m. [27-29]. In non-stressful situations, both CRH and AVP are released in a circadian, pulsatile manner, occurring roughly two to three times per hour [26]. During the early morning hours, the amplitude of CRH and AVP pulses increases, leading to bursts of ACTH and cortisol secretion. High cortisol level may potentially be associated with reduced feedback inhibition of the HPA system, while low cortisol level is linked to increased feedback inhibition [30-32]. Beyond variations in cortisol levels, the changes in the frequency of ultradian oscillations have been observed in certain patients with conditions like depression, post-traumatic stress disorder (PTSD), post-infection fatigue, and chronic fatigue syndrome [33-36]. These variations can result in malfunctions of the HPA axis, making it an area for deeper exploration through dynamic modeling.

Several HPA axis mathematical models explained the dysfunction of the HPA axis dynamics in last decades. In this study, a comprehensive review of existing models and proceeded to replicate and compare five of them using a data set provided by Young et al. [55]. Our analysis reveals that none of these models demonstrated a strong fit to the validation dataset. In response to these limitations, we have extended the minimal mechanistic model of the HPA axis through the integration of the role of the hippocampus in simulating the hypothalamus to produce CRH and a daylight function representing time-related information and circadian rhythms (known as zeitgebers) from the suprachiasmatic nucleus (SCN) of the hypothalamus in the form of trigonometric functions. Furthermore, we introduce a short negative feedback loop from the pituitary to the hypothalamus and an additional signaling pathway involving Arginine Vasopressin (AVP) in the production and secretion of ACTH. This integration allows us to estimate model parameters that helped to a significant reduction in the mean absolute percent error. By sensitivity and correlation analyses, this model shows the parameters that exert the most significant influence on the dynamics of cortisol within the system. Although existing advancements in the HPA model have helped to comprehend the dynamics of the HPA Axis, these models are still unable to incorporate major factors such as the impact of genetic and epigenetics, the role of the amygdala in processing and transmitting stress signals, and other minor impacting factors.

**1.2. Background: Mathematical Models of HPA Axis in the Last Three Decades**

In 1994 Gonzalez-Heydrich et al.[37], studied the dynamics of CRH, ACTH, and cortisol by developing the computer model using ODE. The foundational structure of these models has remained fairly consistent since then. The first enhancement to these models was introduced by Liu et al.[38] in 1999. They expanded the model to include CRH, ACTH, free cortisol, cortisol bound to cortisol-binding globulin (CBG), and cortisol bound to albumin. In 2005, Smiljana Jelic [39] and her team proposed a nonlinear differential equation model to differentiate between plasma cortisol levels in humans under normal, homeostatic conditions. This model aimed to describe the activity of the HPA axis system. In 2007, S. Gupta [40] introduced additional complexities to this model, incorporating glucocorticoid receptor (GR) kinetics. This expanded model sought to elucidate the dysregulation seen in Chronic Fatigue Syndrome and revealed the presence of bi-stability within the system. In 2008, first-time Bairagi et al. [41] introduced the delay differential equation (DDE) based model by including delays between the production of ACTH and cortisol and their action. They had shown, both ultradian and circadian oscillations are generated using this model. But, the circadian oscillations required a pulse generator function that modeled the suprachiasmatic nucleus (SCN) of the hypothalamus. J.J. Walker's [42] model in 2010 elucidated the role of a combination of delay with feed-forward, and feedback loops in the pituitary-adrenal system, providing insights into the ultradian pulsatility in the absence of an ultradian source from a supra-pituitary site. The minimal model of the HPA axis was introduced by Frank Vinther et al. [43] in 2011 by incorporating the most widely accepted factors that influence HPA axis dynamics, keeping in mind the negative feedback loop involving cortisol, CRH, and ACTH in a system of three coupled, nonlinear differential equations. Subsequently, In 2012 Sriram et al. [44] presented a model that produced both circadian and ultradian oscillations without external inputs from SCN, revealing the significance of GR interactions. This model consisted of four ordinary differential equations CRH, ACTH, cortisol, and GR availibilty/binding. On further exploration of the HPA axis's dynamics in 2013, Stefen Hui et al. [45] proposed a non-linear model to create a stressor estimator, offering potential applications in stress-related disease treatment. In 2013, Andersen et al. [46] introduced the role of hippocampal GR/MR interactions in a model to produce oscillations. However, this model failed to generate any oscillations for physiologically reasonable parameter values. In 2015 Niyousha et al. [47] presented a review of existing models by comparing concentration datasets of ACTH and cortisol from 17 healthy individuals and extended the model by adding individual fixed effects and an exogenous circadian function after analysing contradictory results of various feedback mechanisms. In 2015, H. Malek et al. [48], published a simple delayed model of the HPA axis that examines how the hypothalamus-pituitary-adrenal (HPA) axis and inflammatory cytokines interact in cases of acute inflammation during both normal and infectious conditions. In 2017, Bangsgaard and Ottesen [49] developed a model that utilized an optimization procedure to fit experimental data obtained from both healthy individuals and those suffering from depression. Through this approach, the model was capable of comparing the parameters associated with depressed participants and those from a control group of healthy subjects, yielding valuable insights. In 2019, Casey Johnson et al. [50] presented the HPA-axis model, which analyzes the linear and nonlinear stability of stationary solutions by incorporating the glucocorticoid receptor but does not take into account the system response delay. So, taking into account a delay in system response explains the mechanism of the HPA axis self-regulatory activities, the existence of periodic solutions, and analyze stability of solutions with respect to time delay values. To explain the dysregulation of stress hormones, In 2020, Uri Alon et al. [51] presented a minimal model of the HPA axis on the timescale of weeks by incorporating changes in the total functional mass of the HPA axis hormone-secreting glands in the situations of alcohol abuse, anorexia, and postpartum. In 2022, Churilov et. al [52] introduced a system of impulsive time-delay differential equations that incorporates the pulsatile release of adrenocorticotropin (ACTH) by the pituitary gland and a time delay for the release of glucocorticoid(cortisol) hormones from the adrenal gland. A computational module coded in Python called VeVaPy for verification and validation using 10 selected HPA axis models was presented by Christopher Parker et al. [53] in 2022 to assist future researchers in verifying mathematical models. In 2023, Caruso et al. [54] incorporated novel signaling pathways of Arginine Vasopressin to enhance the predictive ability of their model. Moreover, it provides deeper insights into cortisol secretion in the form of ultradian rhythm by including pulsatile release of Adrenocorticotropic Hormone combined with negative feedback into the system from glucocorticoid receptors in the model. The details are different approach to HPA Axis mathematical modelling is summarized in the diagrams Fig 1(A)and 1(B).

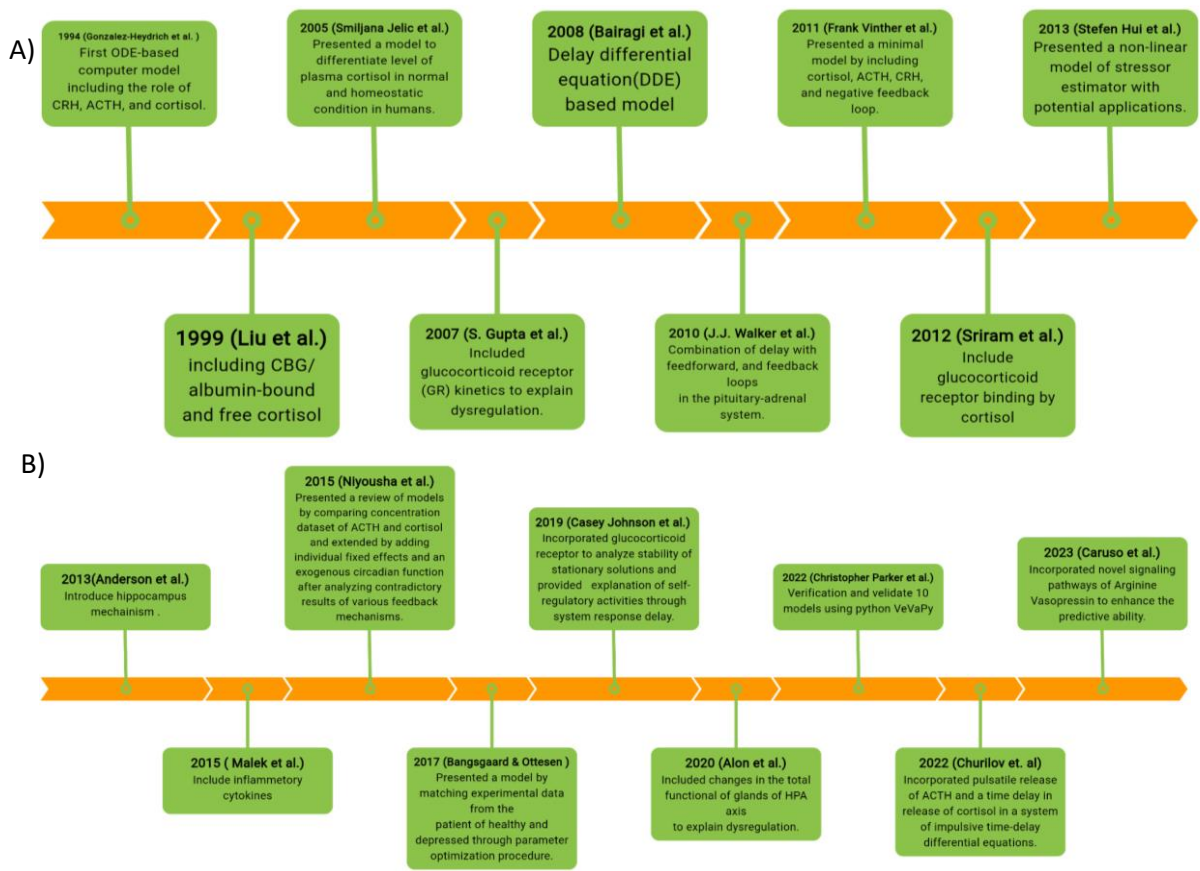

*Figure 1: Timeline of hypothalamus-pituitary-adrenal (HPA) axis ODE based model. Models included, in Chronological order: (A) Gonzalez-Heydrich et al. [37], Liu et al. [38], Smiljana Jelic et al. [39], S. Gupta[40], Bairagi et al. [41], J.J. Walker[42], Frank Vinther et al.[43], Sriram et al. [44], Stefen Hui et al.[45] (B) Andersen et al.[46], Malek et al.[48], Niyousha et al.[47], Bangsgaard and Ottesen [49], Casey Johnson et al.[50], Uri Alon et al.[51], Christopher Parker et al.[53], Churilov et. al [52], Caruso[54].*

**1.3 Model Validation against data:**

Mean Absolute Percentage Error (MAPE) and Root Mean Square Error (RMSE) are essential metrics used for the validation and comparison of predictive models. RMSE focuses on the magnitude of errors in the same units as the data, while MAPE concentrates on the percentage difference between predictions and actual data, making it particularly useful for data with varying ranges and differing units. The unit of hormone concentrations varies across the selected models, specifically between mol/dm3 and ng/ml. Due to this variation in units, MAPE be the most relevant fit measure for analysis. MAPE is expressed as a percentage, with lower values indicating more accurate model predictions. Ideally, a MAPE of 0% would signify perfect predictions, while larger percentages correspond to larger errors.

$$\text{MAPE} = \frac{1}{n}\sum_{i=1}^{n}\left|\frac{Predicted_i - Actual_i}{Actual_i}\right| \times 100$$

RMSE is used to assess the quality of predictive models. RMSE measure in the same units as the data, which allows for easy comparison. A lower RMSE indicates superior model performance and is particularly valuable when the units of measurement are consistent [55].

$$\text{RMSE} = \sqrt{\frac{1}{n}\sum_{i=1}^{n}(Predicted_i - Actual_i)^2}$$

In this study the data of ACTH and Cortisol concentration taken for 24 hours in depressed women obtained by Young et al. [56] is used for validation of five HPA axis models. The validation result is given in table 1. This comparative analysis is a pivotal measure to evaluate the accuracy and reliability of the models.

**Table 1:**

Validation results for Cortisol and ACTH.

| Studies | Cortisol MAPE(%) | ACTH MAPE(%) | Cortisol RMSE | ACTH RMSE |
|---|---|---|---|---|
| Caruso et al.[54] | 59.9236 | 59.0381 | 6.2693 | 12.1951 |
| Andersen et al.[46] | 62.7917 | 47.1386 | 6.3987 | 11.1082 |
| Sriram et al.[44] | 85.1361 | 81.5760 | 6.8293 | 15.5908 |
| Markovic et al.[57] | 72.6983 | 100 | 6.9442 | 17.145 |
| Ben zvi et al.[58] | 94.50 | 99.99 | 8.3161 | 17.1459 |

**2. Model Formulation:**

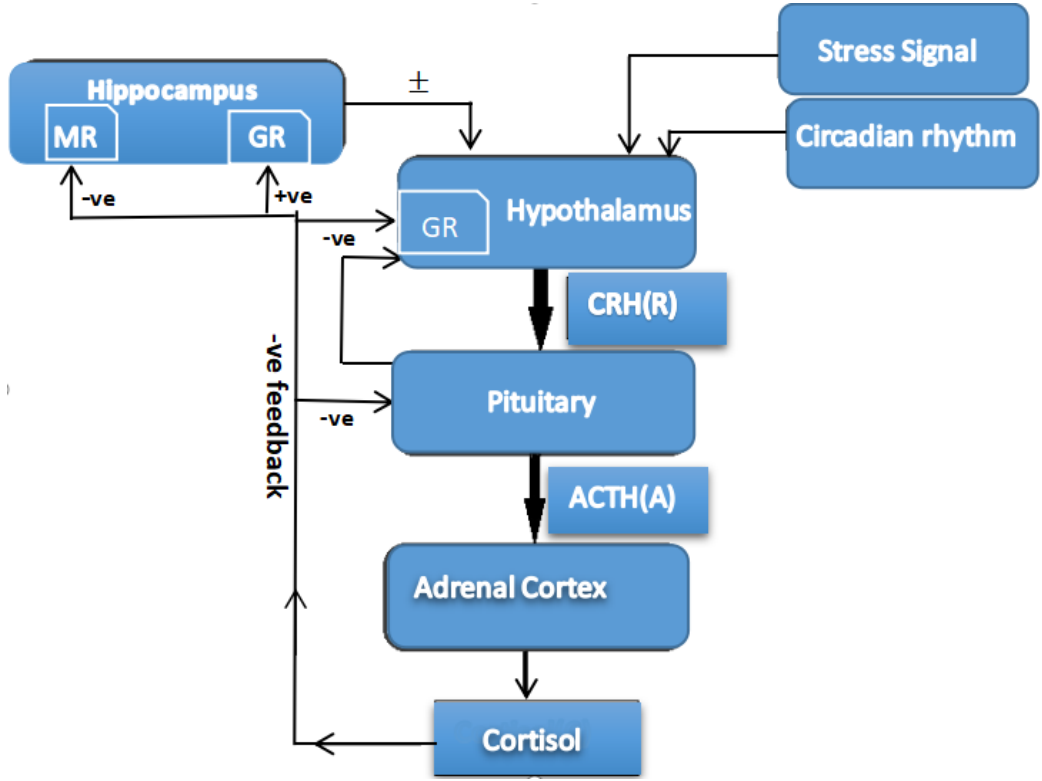

*Figure 2:* Schematic diagram of extended Model

In this study, we explore the intricate mechanisms governing the suprachiasmatic nucleus (SCN) response to zeitgebers, which is deeply intertwined with the day-night cycles and exerts a profound influence on the fluctuations of hormones within the HPA axis. To model the circadian rhythm, we introduce the Day function (D) as described in equation (2.1). Furthermore, this paper introduces a novel addition to the model: the incorporation of the hippocampus mechanism, denoted as $(k_0 = -\psi \frac{C^\delta}{R_C^\delta + C^\delta} + \eta \frac{C^\beta}{R_C^\beta + C^\beta})$. The hippocampus, through a complex interplay involving mineralocorticoid receptors (MR) and glucocorticoid receptors (GR) in conjunction with cortisol, partially controls the release of CRH (Corticotrophin-Releasing Hormone). This augmentation in the model significantly improves its accuracy. Adding the effect of the hippocampus mechanism to the HPA axis model leads to substantial enhancements. Specifically, the Mean Absolute Percentage Error (MAPE) for Cortisol is reduced by 7.84 percentage points (from 59.92% to 51.08%), and for ACTH, it is reduced by

5.88 percentage points (from 59.03% to 53.15%) when $k_4$ changes from 0.0821 to 0.0801 from the original model. The model is composed of a set of interconnected nonlinear differential equations (2.3) that capture the rates of change of CRH (R), ACTH (A), and cortisol (C) within the HPA axis. These equations account for various factors, including positive stimulations, density-dependent stimulations, and negative feedback mechanisms, governing the intricate dynamics of this complex physiological system. The first nonlinear differential equation describes the rate of change of CRH (R) from the hypothalamus, with positive stimulation from the stress signal of ultradian rhythm and circadian rhythm with the daylight function (D). However, it is essential to recognize that the HPA axis operates within a feedback mechanism involving the hypothalamus and pituitary gland. When cortisol levels rise, signals are transmitted to reduce CRH and ACTH production, subsequently curtailing cortisol production when it is no longer required. This intricate feedback loop ensures the system's responsiveness to negative feedback, leading to the inhibition of further CRH release via cortisol from the adrenal gland $\left(1 - \mu \frac{C^\beta}{R_C^\beta + C^\beta}\right)$ and ACTH $\left(1 - \phi \frac{A^\alpha}{R_A^\alpha + A^\alpha}\right)$.

$$D = \frac{1}{11.1}\left[3.9 \sin\left(\frac{\pi t}{720}\right) - \sin\left(\frac{2\pi t}{720}\right) - 1.3 \cos\left(\frac{2\pi t}{720}\right) - 2.8 \cos\left(\frac{\pi t}{720}\right)\right] + 0.4 \quad (2.1)$$

In contrast, the hippocampus also exerts partial control over the release of CRH through a complex interplay of mineralocorticoid receptors (MR) and glucocorticoid receptors (GR) in conjunction with cortisol. Thus, we introduce the hippocampus mechanism to our model $\left(k_0 = -\psi \frac{C^\delta}{R_C^\delta + C^\delta} + \eta \frac{C^\beta}{R_C^\beta + C^\beta}\right)$. CRH (R) is released with a removal rate ($h_1$), contributing to the regulatory dynamics.

$$\frac{dR}{dt} = (k_1 + Dk_2)\left(1 - \phi \frac{A^\alpha}{R_A^\alpha + A^\alpha}\right)\left(1 - \mu \frac{C^\beta}{R_C^\beta + C^\beta} - \psi \frac{C^\delta}{R_C^\delta + C^\delta} + \eta \frac{C^\beta}{R_C^\beta + C^\beta}\right) \quad (2.2a)$$
$$- h_1 R$$

The second equation captures the rate of change of ACTH from the pituitary gland $\left(\frac{dA}{dt}\right)$. This is influenced by constant stimulation from CRH ($k_4$) and arginine vasopressin ($k_3$), combined with density-dependent stimulation $\left(\frac{D^\gamma}{R_D^\gamma + D^\gamma}\right)$, which travels a short distance to the pituitary, leading to ACTH (A) release. On the other hand, negative feedback from cortisol $\left(1 - \rho \frac{C^\beta}{R_C^\beta + C^\beta}\right)$ inhibits further ACTH release. As a result, ACTH (A) is released with a removal rate ($h_2$).

$$\frac{dA}{dt} = \left(k_3 \left(\frac{D^\gamma}{R_D^\gamma + D^\gamma}\right) + k_4 R\right)\left(1 - \rho \frac{C^\beta}{R_C^\beta + C^\beta}\right) - h_2 A \quad (2.2b)$$

The third differential equation assists in measuring the rate of change of cortisol $\left(\frac{dC}{dt}\right)$ from the adrenal gland. ACTH (A) enters the adrenal cortex, providing a positive stimulus ($k_5$) for cortisol release (C) into the bloodstream, with a constant removal rate ($h_3$). This cause to increase in stress levels.

$$\frac{dC}{dt} = k_5 A - h_3 C \quad (2.2c)$$

The set of Equations (2.2) reduces to

$$\frac{dR}{dt} = (k_1 + Dk_2)\left(1 - \phi \frac{A^\alpha}{R_A^\alpha + A^\alpha}\right)\left(1 - \mu \frac{C^\beta}{R_C^\beta + C^\beta} - \psi \frac{C^\delta}{R_C^\delta + C^\delta}\right.$$
$$\left. + \eta \frac{C^\beta}{R_C^\beta + C^\beta}\right) - h_1 R \quad (2.3a)$$

$$\frac{dA}{dt} = \left(k_3\left(\frac{D^\gamma}{R_D^\gamma + D^\gamma}\right) + k_4 R\right)\left(1 - \rho \frac{C^\beta}{R_C^\beta + C^\beta}\right) - h_2 A \qquad (2.3b)$$

$$\frac{dC}{dt} = k_5 A - h_3 C \qquad (2.3c)$$

**Table 2. Details of Parameters values and associated references**

| Variables | Represents | Values | Source and Ref | Units |
|---|---|---|---|---|
| $h_1$ | Removable rate of CRH | 0.1732/min | (Andersen et al, 2013) [46] | µg/dL |
| $h_2$ | Removable rate of ACTH | 0.0315/min | (Andersen et al, 2013) [46] | pg/mL |
| $h_3$ | Removable rate of Cortisol | 0.0105/min | (Andersen et al, 2013) [46] | µg/dL |
| $R_C$ | Cortisol Saturation Constant | 1.12 | Caruso [54] | µg/dL |
| $R_A$ | ACTH Saturation Constant | 0.78 | Caruso [54] | Pg/mL |
| $R_D$ | AVP saturation constant | 1.3 | Caruso [54] | dimensionless |
| α | Hill Coefficient | 4 | Caruso [54] | dimensionless |
| β | Hill Coefficient | 3 | Caruso [54] | dimensionless |
| γ | Hill Coefficient | 3 | Caruso [54] | dimensionless |
| ϕ | Level of inhibition of CRH by ACTH | 0.160 | Caruso [54] | dimensionless |
| ψ | Level of inhibition of CRH by Cortisol | 0.5 | Caruso [54] | dimensionless |
| ξ | $\phi - \mu$ | 2 | Andersen et. al. [46] | dimensionless |
| ρ | Level of inhibition of ACTH by Cortisol | 0.304 | Caruso [54] | dimensionless |

| $k_1$ | Stimulation of CRH from outside of the SCN | 0.5703 | Fitted | pg/mLmin |
| $k_2$ | Stimulation of CRH by Hypothalamus | 0.4342 | Fitted | pg/mLmin |
| $k_3$ | Stimulation of CRH by AVP | 0.2166 | Fitted | pg/mLmin |
| $k_4$ | Stimulation of ACTH by CRH | 0.0821 | Fitted | 1/min |
| $k_5$ | Stimulation of Cortisol by ACTH | 0.00430 | Fitted | $10^5$/min |

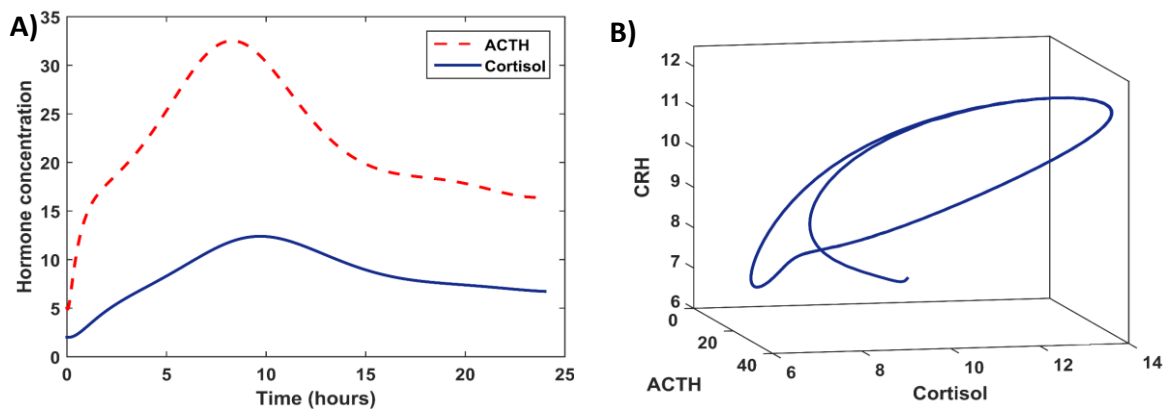

*Figure 3:* Simulation result of ACTH and Cortisol concentration vary with time where 0 hour corresponds to mid night (sub plot A). B) The periodic solution of Cortisol and ACTH production in every 24 hours.

The HPA axis, responsible for regulating ACTH and cortisol levels, exhibits a rhythmic pattern over a 24-hour period. This cyclic behaviour is vital for maintaining the body's homeostasis and adapting to daily challenges. In the morning, there is a notable peak in both ACTH and cortisol levels, often referred to as the "cortisol awakening response." This surge, occurring around *6-8 a.m.*, readies the body for the demands of the day, offering energy and stress management support. As the day unfolds, ACTH and cortisol levels gradually decline.

By late evening and into the night, cortisol levels reach their lowest point, allowing the body to enter a state of rest and repair. During these night time, processes like tissue repair, cell growth, and immune system function are optimized. HPA axis mechanism can be visually elucidated using a time series graph. The x-axis would represent the 24-hour time span, starting at midnight (t=0) and ending at the next midnight. The y-axis would depict the concentration levels of ACTH and cortisol. Peaks in the morning and troughs at night should be clearly visible, illustrating the circadian rhythm.

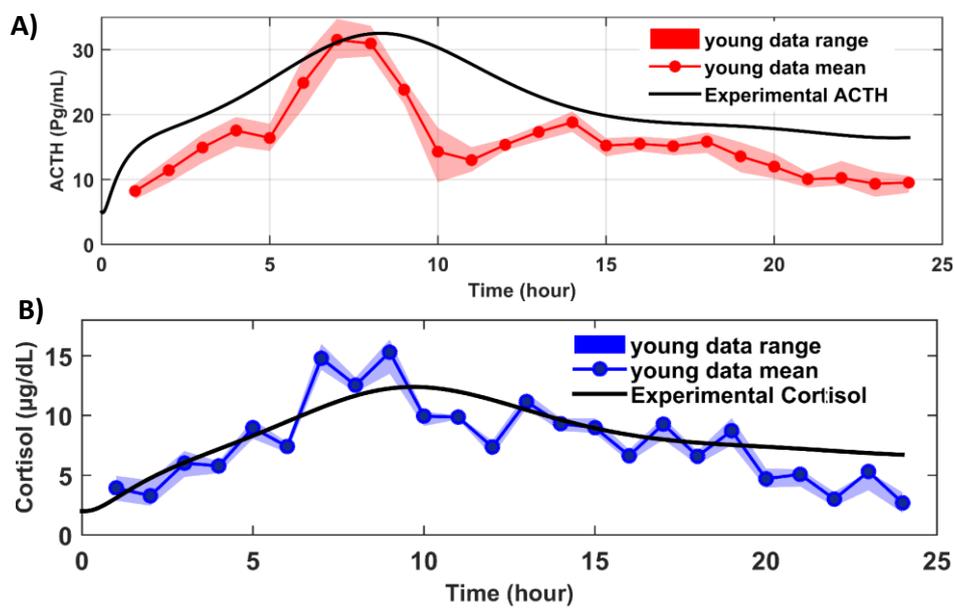

*Figure 4:* *The plot a) shows the Cortisol concentrations predicted by the model displayed in black and for comparison the young et al.[56] data displayed in blue. The plot B) shows the ACTH concentration predicted by the model displayed in black and for comparison the young et al.[56]data displayed in red.*

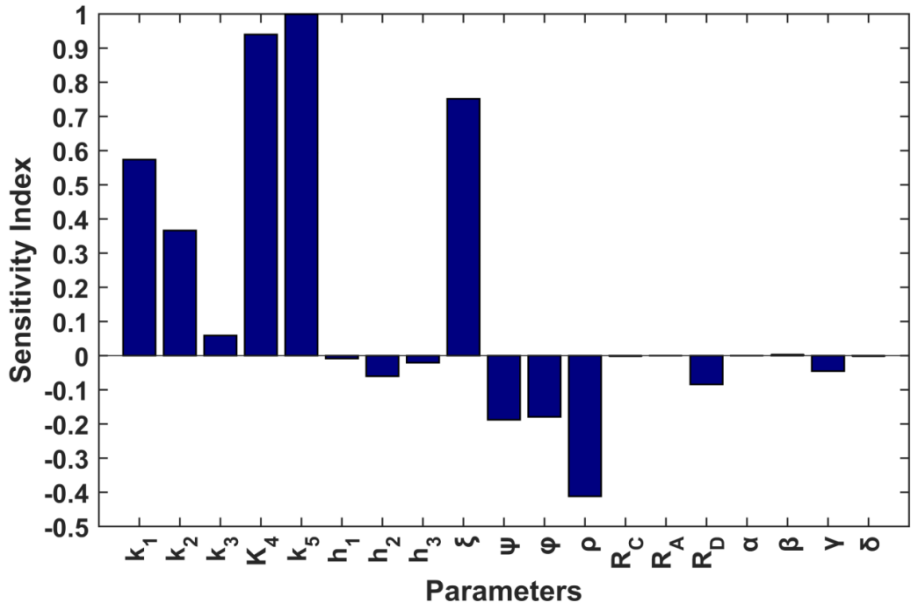

*Figure 5:* Sensitivity indices for optimal set of parameters. Sensitivity analysis shows the sensitivity of cortisol for the optimal set of parameters.

### 3. Sensitivity and Correlation Analysis

To gain deeper insights into the intricate dynamics of cortisol concentration predicted by our computational model, we conducted a comprehensive sensitivity and correlation analysis. This analysis allows us to understand how changes in various parameters affect the concentration of cortisol in the model output. Moreover, it sheds light on potential relationships and dependencies between these parameters. The primary objective of the sensitivity analysis was to evaluate the model's responsiveness to alterations in specific parameters while keeping the rest of the parameters constant. We conducted this analysis using the optimal parameter set as detailed in Table 1. In this investigation we use relative sensitivity indices (SI) to quantify the impact of parameter variations on cortisol concentration. The relative sensitivity index (SI) for the model's cortisol output, as calculated for the optimal parameter set, can be expressed as follows:

$$SI(p_i) = \frac{\partial C}{\partial p_i} \frac{p_i}{C}$$

The sensitivity values for each parameter were determined and subsequently ordered in descending order as $k_5 > k_4 > \xi > k_1 > \rho > k_2 > \psi > \phi > R_D > h_2 > k_3 > \gamma > h_3 > h_1 > \beta > R_C > \delta > R_A > \alpha$. The model exhibited negligible sensitivity to changes in some parameters while displaying a highly responsive behaviour to variations in others. Specifically, we found that the model was almost unresponsive to changes in parameter "$\alpha$" and "$R_A$" but exhibited remarkable sensitivity to changes in parameter "$k_5$" and "$k_4$". This disparity in sensitivity highlights the importance of parameter "$k_5$" and "$k_4$" in shaping the dynamics of cortisol concentration, suggesting that it plays a pivotal role in our model's predictions. Correlation analysis helps us understand how changes in one parameter may be linked to changes in another. A strong positive correlation was observed between the sensitivities of parameters $k_5$ and $k_4$, denoted as $r(k_5,k_4)=+1$. This finding suggests that changes in these two parameters tend to a cooperative effect on cortisol dynamics while negative correlations suggest an inverse relationship.

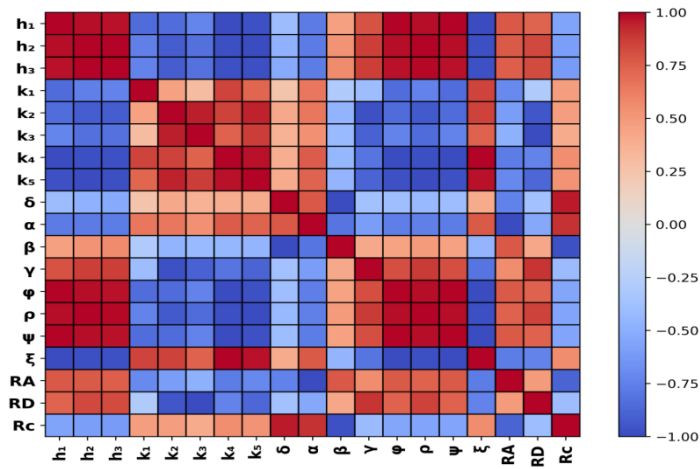

*Figure 6: Correlation matrix. The correlation values between the dynamics sensitivities for the all the parameters are shown with the diagonal being self-correlated. The level of correlation are differently shaded as shown on the horizontal bar on the right that ranges from highly correlated (+1) to anti-correlated (-1).*

**4. Limitations and Future Research Directions**

The base model discussed in this paper does not provide a highly accurate fit for hormone levels between *8 am* and *12 pm*. It fails to capture the quick decline followed by a small second peak observed in the data during this time period. While this limitation does not drastically affect the interpretation of results, it could impact findings related to the generation of ultradian rhythms.

Although the development in mathematical models of the HPA axis has significantly progressed to encompass the various factors that influence stress responses, these models still lack following several major impacting factors. Genetic variation alters hormones release by changing the structure and function of the HPA axis. Epigenetic modifications in the form of DNA methylation and histone acetylation change gene expression within the HPA axis, potentially increasing or decreasing stress reactivity. The function of the amygdala in processing and transmitting stress signals to the hypothalamus play essential role in both the emotional and physiological aspects of stress. Furthermore, the integration of smaller factors such as physiological and psychological individuality, gender-specific hormonal differences, growth hormone, ghrelin, leptin, pituitary-thyroid interactions, the effects of specific medications, and hormonal therapies into the mathematical framework. It could enable us to create a holistic mathematical model capable of capturing the intricate complexities involved in activating the HPA axis during stress responses.

5. **Discussions**

In this study delves into the intricate mechanisms that govern the Suprachiasmatic Nucleus (SCN) response to zeitgebers, shedding light on its profound influence on the Hypothalamus-Pituitary-Adrenal (HPA) axis and its role in regulating hormone fluctuations. We have introduced a novel addition to the model, the hippocampus mechanism, which significantly enhances the accuracy of our model, reducing Mean Absolute Percentage Error (MAPE) for cortisol by 7.84 percentage points and for ACTH by 5.88 percentage points.Our model, composed of a set of interconnected nonlinear differential equations, captures the complex dynamics of the HPA axis, accounting for positive stimulations, density-dependent stimulations, and negative feedback mechanisms. It illustrates the intricate feedback loop between the hypothalamus, pituitary gland, and adrenal gland, which is vital for maintaining physiological balance.The HPA axis exhibits a rhythmic pattern over a 24-hour period, crucial for the body's daily adaptation to challenges. The cortisol awakening response in the morning prepares the body for the day, while cortisol levels decrease in the evening, allowing the body to rest and repair during the night. The sensitivity and correlation analyses have provided valuable insights into the model's responsiveness to parameter variations and potential dependencies between these parameters. Certain parameters, such as the removal rate of CRH and cortisol saturation constants, were found to significantly impact cortisol dynamics. Correlation analysis further revealed cooperative or inverse relationships between specific parameters. Our work represents a significant contribution to the field and paves the way for future research directions. These directions include exploring genetic influences on HPA axis regulation, developing personalized medication approaches, investigating neural circuits within the axis, and elucidating the connections between HPA axis dysregulation and stress-related diseases. The knowledge gained from this research will undoubtedly shape and guide the future of HPA axis research, offering new possibilities for understanding and addressing stress-related health issues.